\documentclass[letterpaper]{amsart}
\usepackage{amssymb}

\newtheorem{thm}{Theorem}
\newtheorem{lem}[thm]{Lemma}

\newtheorem{cor}[thm]{Corollary}

\theoremstyle{definition}

\newtheorem*{rem}{Remark}
\newtheorem{exmp}{Example}

\newtheorem*{conj}{Conjecture}

\newcommand{\bQ}{\mathbb Q}
\newcommand{\bZ}{\mathbb Z}
\newcommand{\bF}{\mathbb F}

\begin{document}
\title[A Valuation Criterion for Normal Bases]{A Valuation Criterion for Normal Bases in
  Elementary Abelian Extensions}

\author[Nigel P. {\sc Byott}]{{\sc Nigel P.} Byott}
\address{Nigel P. {\sc Byott}\\
School of Engineering, Computer Science and Mathematics \\
University of Exeter \\
Exeter EX4 4QE \\
United Kingdom}
\email{N.P.Byott@ex.ac.uk}

\author[G. Griffith {\sc Elder}]{{\sc G. Griffith} Elder}
\address{G. Griffith {\sc Elder}\\ Department of Mathematics \\
University of Nebraska at Omaha\\ Omaha, NE 68182-0243 U.S.A.}
\email{elder@vt.edu} \curraddr{Department of Mathematics, Virginia
Tech, Blacksburg VA 24061-0123}

\thanks{Elder was partially supported by NSF grant DMS-0201080.}
\date{September 26, 2006}

\begin{abstract}
Let $p$ be a prime number and let $K$ be a finite extension of the
field $\bQ_p$ of $p$-adic numbers. Let $N$ be a fully ramified,
elementary abelian extension of $K$. Under a mild hypothesis
on the extension $N/K$, we show that every element of $N$ with valuation
congruent mod $[N:K]$ to the largest lower ramification number of
$N/K$ generates a normal basis for $N$ over $K$. 
\end{abstract}

\bibliographystyle{amsalpha} 

\maketitle

\section{Introduction}
The Normal Basis Theorem states that in a finite Galois extension
$N/K$ there are elements $\alpha\in N$ whose conjugates
$\{\sigma\alpha:\sigma\in \mbox{Gal}(N/K)\}$ provide a vector space
basis for $N$ over $K$. If $K$ is a finite extension of the field $\bQ_p$
of $p$-adic numbers, the valuation $v_N(\alpha)$ of an element $\alpha$ of $N$ is
an important property. We therefore ask whether anything can be said
about the valuation of normal basis generators in this case. We will prove

\begin{thm}
Let $K$ be a finite extension of the $p$-adic numbers, let $N/K$ be a
fully ramified, elementary abelian $p$-extension, and let $b_{\max}$
denote the largest lower ramification number.  If the upper
ramification numbers of $N/K$ are relatively prime to $p$, then every
element $\alpha\in N$ with valuation $v_N(\alpha)\equiv b_{\max}\bmod
[N:K]$ generates a normal field basis. Moreover, no other equivalence class has
this property: given any integer $v$ with $v \not \equiv  b_m
\bmod [N:K]$, there is an element $\rho_v \in N$ with $v_N(\rho_v)=v$
which does not generate a normal basis.
\end{thm}

This result arose out of work on the Galois module structure of ideals
in extensions of $p$-adic fields. For such extensions, it has been
found that the usual ramification invariants are, in general,
insufficient to determine Galois module structure, and thus that there
is a need for a {\em refined ramification filtration}
\cite{elder:onebreak, elder:newbreaks, elder:necbreaks}. This refined
filtration is defined for elementary abelian $p$-extensions and
requires elements that generate normal field bases. Such elements are
provided by Theorem 1. Recent work \cite{elder:onedim} suggests that
what is known for $p$-adic fields should also hold in the analogous
situation in characteristic $p$, where $K$ is a finite extension of 
$\bF_p(X)$. Here $\bF_p$ denotes the finite field with $p$ elements,
and $X$ is an indeterminate. We therefore make the

\begin{conj}
Theorem 1 holds when $K$ is a finite extension of $\bF_p(X)$  as well.
\end{conj}
\section{Preliminary Results}

Let $K$ be a finite extension of the field $\bQ_p$ of $p$-adic
numbers, and let $N/K$ be a fully ramified, elementary abelian
$p$-extension with $G=\mbox{Gal}(N/K)\cong C_p^n$.  Use subscripts to
denote field of reference. So $\pi_N$ denotes a prime element in $N$,
$v_N$ denotes the valuation normalized so that $v_N(\pi_N)=1$, and
$e_K$ denotes the absolute ramification index.  Let $\mbox{Tr}_{N/K}$
denote the trace from $N$ down to $K$.  For each integer $i \geq -1$,
let $G_i=\{\sigma\in G: v_N((\sigma-1)\pi_N)\geq i+1\}$ be the
$i$th ramification group \cite[IV, \S1]{serre:local}. Then $G_{-1}=G_0=G_1=G$,
and the integers $b$
such that $G_b\supsetneq G_{b+1}$ are the lower ramification break (or
jump) numbers. The collection of such numbers, $b_1<\cdots <b_m$, is
the set of lower breaks.  They satisfy $b_1\equiv \cdots \equiv
b_m\bmod p$ \cite[IV, \S2, Prop.~11]{serre:local}, where if $b_m\equiv
0\bmod p$ then the extension $N/K$ is cyclic \cite[IV, \S2, 
Ex.~3]{serre:local}.  Let $g_i=|G_i|$. Then the upper ramification
break numbers $u_1<\cdots <u_m$ are given by 
$u_1 = b_1 g_{b_1}/p^n=b_1$ and
$u_i=(b_1g_{b_1}+(b_2-b_1)g_{b_2}+\cdots + (b_i-b_{i-1})g_{b_i})/p^n$
for $i \geq 2$ \cite[IV, \S3]{serre:local}.  

Now by the Normal Basis Theorem, the set
$$\mathcal{NB}=\left\{\rho\in N:\sum_{\sigma\in G}K\cdot \sigma\rho=N
\right\}$$ of normal basis generators is nonempty.  We desire integers
$v\in \bZ$ such that $\{\rho\in N: v_N(\rho)=v\}\subset
\mathcal{NB}$. And so we are concerned by the following

\begin{exmp}
  Suppose $K$ contains a $p$th root of unity $\zeta$, and let $N=K(x)$ with $x^p-\pi_K=0$.
  Let $\sigma$ generate $\mbox{Gal}(N/K)$. Observe that
  $(\sigma-1)x^{pi}=0$ and $\mbox{Tr}_{N/K}x^i=0$ for $p\nmid i$. So
  for each $i\in\bZ$, we have $v_N(x^i)=i$ and
  $x^i\not\in\mathcal{NB}$.  Here $N/K$ has one ramification break
  $b=pe_K/(p-1)$, which is divisible by $p$.  \cite[IV, \S2,
 Ex.~4]{serre:local}.
\end{exmp}

\begin{rem} Fortunately, these extensions
provide the only obstacle. The restriction in Theorem 1 to elementary
abelian extensions with upper ramification numbers relatively prime to
$p$ is a restriction to those extensions that do not contain a cyclic
subfield such as in Example 1 \cite[IV, \S3 Prop.~14]{serre:local}.
\end{rem}

To prove Theorem 1 we need two results.

\begin{lem}
Let $N/K$ be as above with $b_m\not\equiv 0\bmod p$, and let
$t_G=\sum_{i=1}^m b_i\cdot |G_{b_i}\setminus G_{b_i+1}|$.  If $\rho\in
N$ with $v_N(\rho)\equiv b_m\bmod p^n$, then $v_N(\mbox{\rm
Tr}_{N/K}\rho)=v_N(\rho)+t_G$.  Conversely, given $\alpha\in K$ there
is a $\rho\in N$ with $v_N(\rho)=v_N(\alpha)-t_G\equiv b_m\bmod p^n$
such that $\mbox{\rm Tr}_{N/K}(\rho)=\alpha$.
\end{lem}
\begin{proof}
Use induction.  Consider $n=1$ when
$\mbox{Gal}(N/K)=\langle\sigma\rangle$ is cyclic of degree $p$. There
is only one break $b$, which satisfies $b<pe_K/(p-1)$. Given $\rho\in
N$ with $v_N(\rho)\equiv b\bmod p$, we have 
$\mbox{Tr}_{N/K}\rho\equiv
(\sigma-1)^{p-1}\rho\bmod p\rho$.  Since $(p-1)b<pe_K$,
$v_N(\mbox{Tr}_{N/K}\rho)=v_N(\rho)+(p-1)b$.  And given $\alpha\in K$,
use \cite[V, \S3, Lem.~4]{serre:local} to find $\rho\in N$ with
$v_N(\rho)=v_N(\alpha)-(p-1)b$ and $\mbox{Tr}_{N/K}\rho=\alpha$.

Assume now that the result is true for $n$, and consider $N/K$ to be a
fully ramified abelian extension of degree $p^{n+1}$.  Recall
$g_i=|G_i|$. Let $H$ be a subgroup of $G$ of index $p$ with
$G_{b_2}\subseteq H$.  Let $L=N^H$ and note that $N/L$ satisfies our
induction hypothesis. Moreover the ramification filtration of $H$ is
given by $H_i=G_i\cap H$ \cite[IV, \S1]{serre:local} . So $|H_i|=g_i$
for $i>b_1$. Therefore $t_H=b_m(g_{b_m}-1)+b_{m-1}(g_{b_{m-1}}-g_{b_m})+
\cdots +b_1(p^n-g_{b_2})$. Given $\rho\in N$ with $v_N(\rho)\equiv
b_m\bmod p^{n+1}$, by induction
$v_N(\mbox{Tr}_{N/L}\rho)=v_N(\rho)+t_H$.  By the Hasse-Arf Theorem,
$p^{n+1}\mid g_{b_i}(b_i-b_{i-1})$ for $1\leq i\leq m$.  Thus
$t_H\equiv -b_m+p^nb_1\bmod p^{n+1}$ and
$v_L(\mbox{Tr}_{N/L}\rho)\equiv b_1\bmod p$.  Using \cite[IV, \S1,
  Prop.~3 Cor.]{serre:local}, $b_1$ is the Hilbert break for the 
$C_p$-extension $L/K$.  Applying the case $n=1$, we find
$v_N(\mbox{Tr}_{N/K}\rho)=v_N(\rho)+t_H+p^n(p-1)b_1= v_N(\rho)+t_G$.
The converse statement follows similarly, using $t_H+p^n(p-1)b_1=
t_G$.
\end{proof}

The following generalizes a technical relationship used in the proof
of Lemma 2.
\begin{lem}
Let $N/K$ be a fully ramified, noncyclic, elementary abelian
extension with group $G \cong C_p^n$.  Let $H$ be a subgroup of $G$ of index
$p$, and let $L=N^H$. If 
$b_m$ is the largest lower break of $N/K$, $b$ the only break of
$N/L$, and $\rho$ any element of $N$ with $v_N(\rho)\equiv b_m\bmod
p^n$, then $v_L(\mbox{\rm Tr}_{N/L}\rho)\equiv b\bmod p$.
\end{lem}
\begin{proof}
In the proof of Lemma 2, $H\supseteq G_{b_2}$ so that 
$G_{b_1}H/H\subsetneq G_{b_1+1}H/H$ following 
\cite[IV, \S1, Prop.~3, Cor.]{serre:local}, and the break for $G/H$
was $b_1$. Here we have no 
such luxury and we have to involve the upper numbers in our
considerations, although the argument is really no different. Note
that there is a $k$ such that $G^{u_k+1}H/H\subsetneq
G^{u_{k}}H/H$. Thus $u_k$ is the upper ramification number of
$G/H$. Since there is only one break in the filtration of $G/H$, the
lower and upper numbers for $G/H$ are the same, $b=u_k$.

The ramification filtration for $H$ is given by taking intersections:
$H_j=G_j\cap H$.  Note that $[G_{b_i}:G_{b_i}\cap H]=p$ for $i\leq k$
and $G_{b_i}\subseteq H$ for $i>k$.  Let $h_j=|H_j|$.  Then
$h_j=g_j/p$ for $j\leq b_k$, and $h_j=g_j$ for $j> b_k$.
Now let
$v_N(\rho)=b_m+p^nt$. Following the proof of Lemma 2 and using the
Hasse-Arf Theorem,
\begin{multline*}
v_N(\mbox{\rm Tr}_{N/L}\rho)=b_m+p^nt+b_m(h_{b_m}-1)+b_{m-1}(h_{b_{m-1}}-h_{b_m})+
\cdots +b_1(h_{b_1}-h_{b_2})\\
=p^nt+(b_m-b_{m-1})h_{b_m}+(b_{m-1}-b_{m-2})h_{b_{m-1}}+ \cdots
+(b_2-b_1)h_{b_{2}}+ b_1h_{b_1}\\ \equiv (b_k-b_{k-1})h_{b_k}+ \cdots
+ (b_2-b_1)h_{b_{2}} +b_1h_{b_1} 
\equiv p^nu_k/p
\equiv p^{n-1}b\bmod p^n
\end{multline*}
Therefore $v_L(\mbox{\rm Tr}_{N/L}\rho)\equiv b\bmod p$.
\end{proof}
\section{Main Result}

\begin{proof}[Proof of Theorem 1]
There are two statements to prove. We begin with the first: We assume
the upper breaks satisfy $p\nmid u_i$, and prove that for $\rho\in N$
$$v_N(\rho)\equiv b_m\bmod p^n\implies \rho\in\mathcal{NB}.$$ The
argument breaks up into two cases: the Kummer case where $\zeta\in K$
and the non-Kummer case where $\zeta\not\in K$. Here $\zeta$ is a
nontrivial $p$th root of unity.

We begin with the Kummer case, and start with $n=1$. Let $\sigma$
generate the Galois group, and denote the one ramification number by
$b$.  Since in this case $b$ is also the upper number, $p\nmid
b$. Therefore $\{v_N((\sigma-1)^i\rho:0\leq i <p\}$ is a complete set
of residues modulo $p$. And since $N/K$ is fully ramified, $\rho$
generates a normal basis.  Now let $n\geq 2$ and note that $N=K(x_1,
x_2, \ldots, x_n)$ with each $x_i^p\in K$.  It suffices to show that
$K[G]\rho$ contains each element $y=x_1^{j_1}x_2^{j_2}\cdots
x_n^{j_n}$ with $0\leq j_i\leq p-1$.  For $y=1$ this is clear, since
$\mbox{Tr}_{N/K}(\rho)\in K$. For any other $y$, let $L=K(y)$ and let
$b$ denote the ramification number of $L/K$.  By Lemma 3,
$v_L(\mbox{Tr}_{N/L}(\rho))\equiv b\bmod p$.  Since $b$ is an upper
number of the ramification filtration of $G$, $p\nmid b$. Now apply the
$n=1$ argument, using $\mbox{Tr}_{N/L}(\rho)$ in $L/K$. Thus $y\in
K[G]\rho$.

We now turn to the non-Kummer case with $\zeta\not\in K$. Let
$E=K(\zeta)$, let $E/K$ have ramification index $e_{E/K}$, and let
$F=N(\zeta)$. Then $F/E$ is a fully ramified Kummer extension of
degree $p^n$. Applying Herbrand's Theorem \cite[IV, \S3,
Lem.~5]{serre:local} to the quotient $G={\mbox{Gal}}(N/K)$ of
$\mbox{Gal}(F/K)$, we find that the maximal ramification break of
$F/E$ is $e_{E/K}b_m \not \equiv 0 \bmod p$. The above discussion for the
Kummer case therefore applies to $F/E$.
Suppose now for a contradiction that $\rho\in N$ with $v_N(\rho)\equiv
b_m\bmod p^n$, and that $K[G]\rho$ is a proper subspace of $N$. Then
by extending scalars (noticing that $E$ and $N$ are linearly disjoint
as their degrees are coprime) we have that $E[G]\rho$ is a proper subspace
of $F$. Moreover $v_F(\rho)\equiv e_{E/K}b_m\bmod p^n$. This
contradicts the result already shown for the Kummer extension $F/K$,
completing the proof of the first statement of the theorem.

Consider the second statement: Given any integer $v$ with $v\not\equiv
b_m\bmod p^n$ there is a $\rho_v\in N$ with $v_N(\rho_v)=v$ such that
$\mbox{\rm Tr}_{N/K}\rho_v=0$ and thus $\rho_v\not\in\mathcal{NB}$.

To prove this statement note that given $v\in \bZ$, there is an $0\leq
a_v <p^n$ such that $v\equiv a_vb_m\bmod p^n$, since $p\nmid b_m$. If
$a_v\neq 1$ we will construct an element $\rho_v\in N$ with
$v_N(\rho_v)=v$ and $\mbox{Tr}_{N/K}\rho_v=0$. To begin, observe that
there is a integer $k$ such that $0\leq k\leq n-1$, $a_v\equiv 1\bmod
p^k$ and $a_v\not\equiv 1\bmod p^{k+1}$.  Recall $g_i=|G_i|$.  Since
the ramification groups are $p$-groups with $g_{i+1}\leq g_i$, there
is a Hilbert break $b_s$ such that $g_{b_s+1}<p^{k+1}\leq g_{b_s}$.
For $i=k, k+1$ choose $H_i$ with $|H_i|=p^i$ and $G_{b_s+1}\subset
H_k\subset H_{k+1} \subseteq G_{b_s}$.  Recall from Lemma 2 the
expression for $t_G$, and note that $t_{H_k}=
b_m(g_{b_m}-1)+b_{m-1}(g_{b_{m-1}}-g_{b_m})+ \cdots
+b_s(p^k-g_{b_{s+1}})\equiv -b_m+b_sp^k\bmod p^n$.  Let $L=N^{H_k}$.
Since $a_v\not\equiv 1\bmod p^{k+1}$, $a_v\equiv 1+rp^k\bmod p^{k+1}$
for some $1\leq r\leq p-1$.  Using the fact that $b_s\equiv b_m\bmod
p$, $a_vb_m+t_{H_k}\equiv (r+1)b_mp^k\bmod p^{k+1}$.  Since $p^k\mid
v_N(\alpha)$ for $\alpha\in L$, we can choose $\alpha\in L$ with
$v_N(\alpha)=v+t_{H_k}-rp^kb_s$. So $v_L(\alpha)\equiv b_s\bmod p$.
Let $\sigma\in G$ so that $\sigma H_k$ generates $H_{k+1}/H_k$.
Therefore $v_N((\sigma-1)^r\alpha)=v+t_{H_k}$. Now using Lemma 2, we
choose $\rho_v\in N$ such that $v_N(\rho_v)=v$ and
$\mbox{Tr}_{N/L}\rho_v=(\sigma-1)^r\alpha$.  Since
$(1+\sigma+\ldots+\sigma^{p-1})\mbox{Tr}_{N/L}\rho_v=0$, we have
$\mbox{Tr}_{N/K}\rho_v=0$.
\end{proof}

\begin{cor}
Let $N/K$ be a fully ramified, elementary abelian extension of degree
$p^n$ with $n>1$ and one ramification break, at $b$.  If $\rho\in N$ with
$v_N(\rho)\equiv b\bmod p^n$, then $\rho\in \mathcal{NB}$.
\end{cor}

\bibliography{bib}
\end{document}